\documentclass[a4paper,12pt]{article}
\usepackage[utf8x]{inputenc}
\usepackage{arabtex}
\usepackage{hyperref}
\usepackage{afterpage}
\usepackage{xcolor}
\usepackage{relsize}
\usepackage{mathrsfs}
\usepackage[all]{xy}
\usepackage{amsmath,amssymb,amsfonts,amscd, graphicx, epstopdf, latexsym, verbatim, multirow, color}
\usepackage{epsfig}
\usepackage{geometry} 
\usepackage{adjustbox}

\newcommand{\F}{\mathcal F}

\newcommand{\Q}{\mathbf Q}
\newcommand{\Z}{\mathbf Z}

\def\xor{\underline{\lor}}
 \renewcommand{\R}{\mathbf R}

\newcommand{\psl}{\mathrm{PSL}_2(\Z)}
\newcommand{\pgl}{\mathrm{PGL}_2(\Z)}

\newtheorem{theorem}{Theorem}
\newtheorem{notatdef}[theorem]{DEFINITION-NOTATION}
\newtheorem{corollary}[theorem]{Corollary}
\newtheorem{lemma}[theorem]{Lemma}
\newtheorem{proposition}[theorem]{Proposition}

\usepackage{color} 

\newcommand\nt[1]{\textcolor{red}{{ #1}}}

\newcommand{\Jimm}{\mathbf J}

\renewcommand{\Delta}{{\rm Discr}}

\newcommand{\sherh}[1]{\fboxsep=0pt\setlength{\fboxrule}{1pt}
\begin{center}
   \fbox{\colorbox{green}{
         \begin{minipage}[t]{13cm}
            #1
         \end{minipage}
      }
   }
\end{center}}
\newcommand{\sherhh}[1]{\fboxsep=0pt\setlength{\fboxrule}{1pt}
\begin{center}
   \fbox{\colorbox{yellow}{
         \begin{minipage}[t]{13cm}
            #1
         \end{minipage}
      }
   }
\end{center}}

\newcommand{\sherhhh}[1]{\fboxsep=0pt\setlength{\fboxrule}{1pt}
\begin{center}
   \fbox{\colorbox{red}{
         \begin{minipage}[t]{13cm}
            #1
         \end{minipage}
      }
   }
\end{center}}
\newcommand{\unut}[1]{#1}

\renewcommand{\sherh}[1]{}\renewcommand{\sherhh}[1]{}\renewcommand{\sherhhh}[1]{}\renewcommand{\nt}[1]{}\renewcommand{\unut}[1]{}
\begin{document}

\title{On the involution of the real line induced by Dyer's outer automorphism of PGL(2,Z)}
\author{A. Muhammed Uluda\u{g}$^*$, Hakan Ayral\footnote{
{Galatasaray University, Department of Mathematics,}
{\c{C}{\i}ra\u{g}an Cad. No. 36, 34357 Be\c{s}ikta\c{s}}
{\.{I}stanbul, Turkey}}}

\maketitle

\centerline{ \it Dedicated to Yılmaz Akyıldız, who shared our enthusiasm about jimm}

\begin{abstract}
We study the  involution of the real line induced by the outer automorphism of the extended modular group PGL(2,Z). 
This `{modular}' involution is discontinuous at rationals but satisfies a surprising collection of functional equations. It preserves the set of real quadratic irrationals mapping them in a non-obvious way to each other. It commutes with the Galois action on real quadratic irrationals.

More generally, it preserves set-wise the orbits of the modular group, thereby inducing an involution of the moduli space of real rank-two lattices.  We give a description of this involution as the boundary action of a certain automorphism of the infinite trivalent tree.  It is conjectured that algebraic numbers of degree at least three are mapped to transcendental numbers under this involution.
\end{abstract}

\section{Introduction}\label{sec:introduction}
It is (it seems not very well-) known that  the group $\pgl$ has an involutive outer automorphism, which was discovered by Dyer in the late 70's \cite{dyer1978automorphism}. 
It would be very strange if this automorphism had no manifestations in myriad contexts where $\pgl$ or its subgroups play a major role. 
Our aim in this paper is to elucidate one of these manifestations.

This paper contains some of the results from our arxiv preprint \cite{jimmarxiv} which could not be published due to its length and notational issues. For some complementary results, see \cite{jimmlebesgue}.

Let $\widehat{\R}:=\R\cup \{\infty\}$. The manifestation in question of $\Jimm$ is a map $\Jimm_\R: \widehat{\R}\rightarrow \widehat{\R}$.
Denoting the continued fractions in the usual way
$$
[n_0,n_1,\dots]=n_0+\cfrac{1}{n_1+\cfrac{1}{{\dots}}},
$$
one has, for an irrational number $[n_0,n_1,\dots]$ with $n_0, n_1, \dots \geq 2$,
\begin{equation}
\Jimm_\R([n_0,n_1,n_2,\dots])=[1_{n_0-1},2,1_{n_1-2},2,1_{n_2-2},\dots]
\end{equation}
where $1_k$ is the sequence ${1,1,\dots, 1}$ of length $k$.
This formula remains valid for $n_0, n_1, \dots \geq 1$, if the emerging $1_{-1}$'s are eliminated  in accordance with the rule $[\dots m, 1_{-1},n,\dots]=[\dots m+n-1,\dots]$ and $1_0$ with the rule 
$[\dots m, 1_{0},n,\dots]=[\dots m,n,\dots]$. See page \pageref{examples} below for some examples.

It is possible to extend this definition of $\Jimm_\R$ to all of $\widehat{\R}$. If we ignore rationals and the noble numbers (i.e. numbers in the $\pgl$-orbit of the golden section, $\Phi:=(1+\sqrt{5})/{2}$), then $\Jimm_\R$ is an involution.
It is well-defined and continuous at irrationals,  but two-valued and discontinuous at rationals (Theorem~\ref{jimmcontini}).  

\medskip\noindent
{\bf Guide for notation}. In what follows,
$\Jimm$ denotes Dyers' outer automorphism of $\pgl$,
$\Jimm_\F$ denotes the automorphism of the tree $|\F|$ induced by $\Jimm$,
$\Jimm_{\partial\F}$ denotes the homeomorphism of the boundary $\partial\F$ induced by $\Jimm_\F$,  and
$\Jimm_\R$ denotes the involution of $\R$ induced by  $\Jimm_{\partial\F}$.
However, we reserve the right to drop the subscript and simply write $\Jimm$ when we think that confusion won't arise.

\subsubsection*{Functional equations} 
The involution $\Jimm_\R$ shares the privileged status of the three involutions generating the extended modular group 
$\pgl$,
$$
K:x\rightarrow 1-x, \quad U:x\rightarrow 1/x, \quad V:x\rightarrow -x,
$$ 
as it interacts in a very harmonious way with them. 
Indeed, removing the rationals and noble numbers from its domain, $\Jimm_\R$ satisfies the following functional equations:

\begin{center}
\fbox{\begin{minipage}{13cm}

\vspace{2mm}
$$
\mbox{(FE:0)} \quad \Jimm(\Jimm(x))=x
$$
$$
\mbox{(FE:I)}\quad \Jimm U=U\Jimm \iff 
\Jimm\bigl(\frac{1}{x}\bigr)=\frac{1}{\Jimm(x)}
$$
$$
\mbox{(FE:II)} \quad\Jimm V= UV\Jimm\iff 
\Jimm(-x)=-\frac{1}{\Jimm(x)}
$$
$$
\mbox{(FE:III)} \quad \Jimm K=K\Jimm
\iff \Jimm(1-x)=1-\Jimm(x)
$$

\vspace{1mm}
\end{minipage}}
\end{center}

The following set of functional equations are derived from the above ones:

\begin{center}
{\fbox{\begin{minipage}{13cm}

$$
\mbox{(FE:IV)} \quad \Jimm KV=KUV \Jimm
\iff \Jimm(1+x)=1+\frac{1}{\Jimm(x)}
$$

\vspace{2mm}
$$
\mbox{(FE:V)} \quad 
(\Jimm_\R M \Jimm_\R) (x)=(\Jimm M)(x) \quad \forall M\in \pgl$$

\vspace{2mm}
$$
\mbox{(FE:VI)} \quad 
\Jimm_\R (M x)=\Jimm M \Jimm_\R(x) \quad \forall M\in \pgl
$$

\vspace{1mm}
\end{minipage}}}
\end{center}

\bigskip
 (FE:I) and (FE:III) states that $\Jimm$ is covariant  with the operators $U$ and $K$.
 (FE:IV) is deduced  from (FE:II) and (FE:III). The most general form of the equations is (FE:VI) and says that $\Jimm_\R$ conjugates the M\"obius map $M\in \pgl$ to 
the M\"obius map $\Jimm(M)\in \pgl$. It is deduced from (FE:I-III) by using the involutivity of $\Jimm$. 

Before attempting to play with them, beware that the equations are not consistent on $\Q$. 
On the set of noble numbers $\Jimm$ is 2-to-1 and is not involutive. 

$\Jimm$ preserves the ``{real-multiplication locus}", i.e. the set of real quadratic irrationals. It does so in a non-trivial way, though it preserves setwise the $\pgl$-orbits of real quadratic irrationals. More generally $\Jimm$ preserves setwise the $\pgl$-orbits on $\widehat{\R}$ thereby inducing an involution of the moduli space of real rank-2 lattices, $\widehat{\R}/\pgl$. For a description of its fixed points, see Proposition~\ref{fixed}.

One can use (FE:I)-(FE:III) to directly define and study $\Jimm_\R$. However, the nature of $\Jimm_\R$ is best understood by considering it as a homeomorphism of the boundary of the Farey tree, induced by an automorphism of the tree.  See the coming sections for details.

Finally, the following consequence of the functional equations is noteworthy:
$$
{\frac{1}{x}+\frac{1}{y}=1\iff \frac{1}{\Jimm(x)}+\frac{1}{\Jimm(y)}=1}.
$$
(See page \pageref{twovariable} for a complete set of two-variable functional equations.) 
Hence $\Jimm$ sends harmonic pairs of numbers to harmonic pairs, 
inducing a duality of Beatty partitions of positive integers. See \cite{jimmarxiv} for some details.
\subsubsection*{Some examples.} \label{examples}
Here is a list of assorted values of $\Jimm$. Recall that $\Phi$ denotes the golden section.
\begin{align}\label{jimmphi}
\Phi=[\overline{1}]\implies\Jimm(\Phi)=\infty=\Jimm\bigl(-{1}/{\Phi}\bigr),
\end{align}
where by $\overline{v}$ we denote the infinite sequence $v,v,v, \dots$, for any  finite sequence  $v$. 
From \ref{jimmphi} by using (FE:II) we find 
\begin{align}\label{jimmminusphi}
\Jimm(-\Phi)=\Jimm\bigl({1}/{\Phi}\bigr)=0.
\end{align}
Repeated application of (FE:IV) gives
\begin{align}\label{jimmenplusphi}
\Jimm(n+\Phi)=\Jimm([n+1,\overline{1}])={F_{n+1}}/{F_{n}}
\end{align}
where $F_n$ denotes the $n$th Fibonacci number. One has
$$
\Jimm([1_n,2,\overline{1}])=[n+1,\infty]=n+1
$$
For the number $\sqrt{2}$ we have something that looks simple
$$
\sqrt{2}=[1,\overline{2}]\implies\Jimm(\sqrt{2})=[\overline{2}]=1+\sqrt{2}
$$
but this is not typical as the next example illustrates:
$$
\Jimm({(3+5\sqrt{2})}/{7})=
\Jimm([1,\overline{2, 3, 1, 1, 2, 1, 1, 1}])=
[2,\overline{2,1,4,5}]
=\frac{-3+2\sqrt{95}}{7} 
$$
In a similar vein, consider the examples
$$
\sqrt{11}=[3;\overline{3, 6}] 
\implies \Jimm(\sqrt{11})=[\overline{1,1,2,1,2,1,1}]=
\frac{15+\sqrt{901}}{26},
$$
$$
\Jimm(-\sqrt{11})=-1/\Jimm(\sqrt{11})=-\frac{26}{15+\sqrt{901}}=\frac{15-\sqrt{901}}{26}.
$$
The last example hints at the following result
\begin{theorem}
The involution $\Jimm$ commutes with the conjugation of real quadratic irrationals; i.e. for every real quadratic irrational $\alpha$ one has
$
\Jimm({\alpha}^*)={\Jimm(\alpha)^*}.
$
\end{theorem}
{\it Proof.} Every such $\alpha$ is the fixed point of some $M\in \psl$, the Galois conjugate $\alpha^*$  being the other root of the equation $Mx=x$. But then $\Jimm(M\alpha)=\Jimm(M)\Jimm(\alpha)=\Jimm(\alpha)$, i.e. $\Jimm(\alpha)$ is a fixed point of  $\Jimm(M)$, the other root being $\Jimm(\alpha)^*$. Finally
$M\alpha=\alpha\implies$ $M\alpha^*=\alpha^*\implies$ $ \Jimm(\alpha^*)=\Jimm(M\alpha^*)=\Jimm(M)\Jimm(\alpha^*)$,
so $\Jimm(\alpha^*)$ is also a fixed point of $\Jimm(M)$, i.e. it must coincide with $\Jimm(\alpha)^*$. \hfill $\Box$

\medskip
To finish, let us give the $\Jimm$-transform of a  non-quadratic algebraic number, 
\begin{eqnarray*}
\Jimm(^3\!\sqrt{2})=
\Jimm([1; 3, 1, 5, 1, 1, 4, 1, 1, 8, 1, 14, 1, 10, 2, 1, 4, 12, 2, 3, 2, 1, 3, 4, 1, \dots])\\
=[2,1,3,1,1,1,4,1,1,4,1_6,3,1_{12},3,1_8,2,3,1,1,2,1_{10},2,2,1,2,3,1,2,1,1,3,\dots ]
\\
=2.784731558662723\dots
\end{eqnarray*}
and the transforms of two familiar transcendental numbers:
\begin{eqnarray*}
\Jimm(\pi)=
\Jimm([3, 7, 15, 1, 292, 1, 1, 1, 2, 1, 3,  \dots])=
[1_2, 2, 1_5, 2, 1_{13}, 3, 1_{290}, 5, 3,  \dots]\\
=1.7237707925480276079699326494931025145558144289232\dots
\end{eqnarray*}
\begin{eqnarray*}
\Jimm(e)=
\Jimm([2,1,2,1,1,4,1,1,6,1,1,8,\dots])= 
[1,3,4,1,1,4,1,1,1,1,\dots, \overline{4,1_{2n}}]\\
=1.3105752928466255215822495496939143349712038085627\dots
\end{eqnarray*}
We have been unable to relate these numbers to other numbers of mathematics.

We made some numerical experiments  on the algebraicity of a few numbers $\Jimm(x)$ where $x$ is an algebraic number of degree $>2$, with a special emphasis on $\Jimm(^3\!\sqrt{2})$.
It is very likely that these numbers are transcendental.

\section{Modular group and its automorphism group}
The {\it modular group} is the projective group $\psl$ of two by two unimodular integral matrices. It acts on the upper half plane by M\"obius transformations. It is the free product of its subgroups generated by $S(z)=-1/z$ and $L(z)= (z-1)/z$, respectively of orders 2 and 3. Thus
$
\psl =\langle S, L \,|\, S^2=L^3=1\rangle.
$
The projective group $\pgl$ consists of integral matrices of determinant $\pm 1$.

\medskip
The fact that $Out(\pgl)\simeq \Z/2\Z$ has a story with a twist:
Hua and Reiner \cite{hua1952automorphisms} claimed in 1952 that $\pgl$ has no outer automorphisms. The error was corrected by Dyer \cite{dyer1978automorphism} in 1978. This automorphism also appears in the work of Djokovic and Miller \cite{djokovic} from about the same time. Dyer also proved that the automorphism tower of $\pgl$ stops here. Note that $\pgl\simeq Aut(\psl)$.

\medskip
Below is a list of some presentations of $\pgl$ and $\Jimm$.
 $$
\begin{array}{|l|l|} \hline &\\[-1ex]
\mbox{\bf Presentation of $\pgl$}\quad &\mbox{\bf The automorphism $\Jimm$}\quad  \\[1ex]\hline &\\[-.9ex]
\langle V,U,K\, |\, V^2=U^2=K^2=(VU)^2=(KU)^3=1 \rangle &(V,U,K) \rightarrow (UV,U,K)\\[1ex]\hline &\\[-.9ex]
\langle V,U,L\, |\, V^2=U^2=(LU)^2=(VU)^2=L^3=1 \rangle&(V,U,L) \rightarrow (UV, U, L)\\[1ex]\hline &\\[-.9ex]
\langle S,V,K\, |\, V^2=S^2=K^2=(SV)^2=(KSV)^3=1 \rangle\quad&(S,V,K) \rightarrow (V,S,K)\\[1ex]\hline
\end{array}
$$

\medskip

It is easy to find the $\Jimm$ of an element of $\pgl$ given as a word in one of the listed presentations. On the other hand, there seems to be no algorithm to compute the $\Jimm$ of an element of $\pgl$ given in the matrix form, other then actually expressing the matrix in terms of one of the presentations above, finding the $\Jimm$-transform, and then computing the matrix. 

The translation $2+x\in \psl$ is sent to  $(2x+1)/(x+1)\in \psl$ under $\Jimm$, 
i.e. $\Jimm$ may send a parabolic element  to a hyperbolic element.

\section{The Farey tree and its boundary} 
\paragraph{The construction of $\F$.}  
From $\psl$ we construct the {\it bipartite Farey tree $\F$} 
tree on which $\psl$ acts, as follows.
The set of edges is $E(\F)=\psl$. Vertices are the left cosets of the subgroups $\langle S\rangle$ and $\langle L\rangle$.  Two distinct vertices $v$ and $v'$ are joined by an edge if and only if $v \cap v'\neq \emptyset$  and in this case the edge between them is the only element in the intersection. Cosets of $\langle S\rangle$ are 2-valent vertices and the cosets of 
$\langle L\rangle$ are 3-valent. Since distinct cosets are disjoint, $\F$ is a bipartite graph. It is connected and circuit-free since $\psl$ is freely generated by $S$ and $L$.

The edges incident to the vertex $\{W,WL,WL^{2}\}$ are $ W, \, WL$ and $ WL^{2} $. These edges inherit a natural cyclic ordering which we fix for all vertices as $(W, WL, WL^2)$. This endows $\F$ with the structure of a ribbon graph. Hence $\F$ is an infinite bipartite tree with a ribbon structure. $M\in \psl$ acts on $\F$ from the left by ribbon graph automorphisms 
by sending the edge labeled $W$ to the edge labeled  $MW$. This action is free on $E(\F)$ but not free on the set of vertices.

\paragraph{The boundary of $\F$.} 
\label{sec:boundary/of/F}
A {\it path} on $\mathcal F$ is a sequence of distinct edges $e_1,e_2,\dots,e_k$  such that $e_i$ and $e_{i+1}$ meet at a vertex, for each $1\leq i< k-1$.  Since the edges of $\F$ are labeled by reduced words in the letters $L$ and $S$, a path in $\F$ is a sequence of reduced words $(W_i)$ in $L$ and $S$, such that $W_i^{-1}W_{i+1}\in \{L, L^2, S\}$ for every $i$.

An {\it end} of $\F$ is an equivalence class of infinite (but not bi-infinite) paths in $\F$, where eventually coinciding paths are considered as equivalent. The set of ends of $\F$ is denoted $\partial \F$. 
The action of $\psl$ on $\F$ extends to  $\partial \F$, the element $M\in \psl$ sending the path $(W_i)$ to the path $(MW_i)$. 

Given an edge $e$ of $\F$ and an end $b$ of $\F$, there is a unique path in the class $b$ which starts at $e$. Hence for any edge $e$, we may  identify the set $\partial \F$ with the set of infinite paths that start at $e$. We denote this latter set by $\partial\F_e$  and endow it with the  topology generated by the {\it Farey intervals} ${\mathcal O}_{e'}$, i.e. the set of paths  through $e'$. The ribbon structure of $\F$ endows $\partial\F_e$ with the structure of a cyclically ordered topological space.  Given  a second edge $e'$ of $\F$, the spaces $\partial\F_e$ and $\partial\F_{e'}$ are  homeomorphic under the order-preserving map which pre-composes with the unique path joining $e$ to $e'$. 
The action of the modular group on the set $\partial\F$ induces an action  of $\psl$ by order-preserving homeomorphisms of $\partial\F_e$, for any choice of a base edge $e$.

\paragraph{The continued fraction map.} 
$\partial\F_e$ is homeomorphic to the Cantor set.
By exploiting its cyclic order structure, we ``smash its holes" as follows. Define a {\it rational end} of $\F$ to be an eventually left-turn or eventually right-turn path. Now introduce the equivalence relation $\sim $ on $\partial\F$ as: left- and right- rational paths which bifurcate from the same vertex are equivalent. On $\partial\F_e$ this equivalence sets equal those points which are not separated by a third point.

\begin{center}
\noindent{\includegraphics[scale=.8]{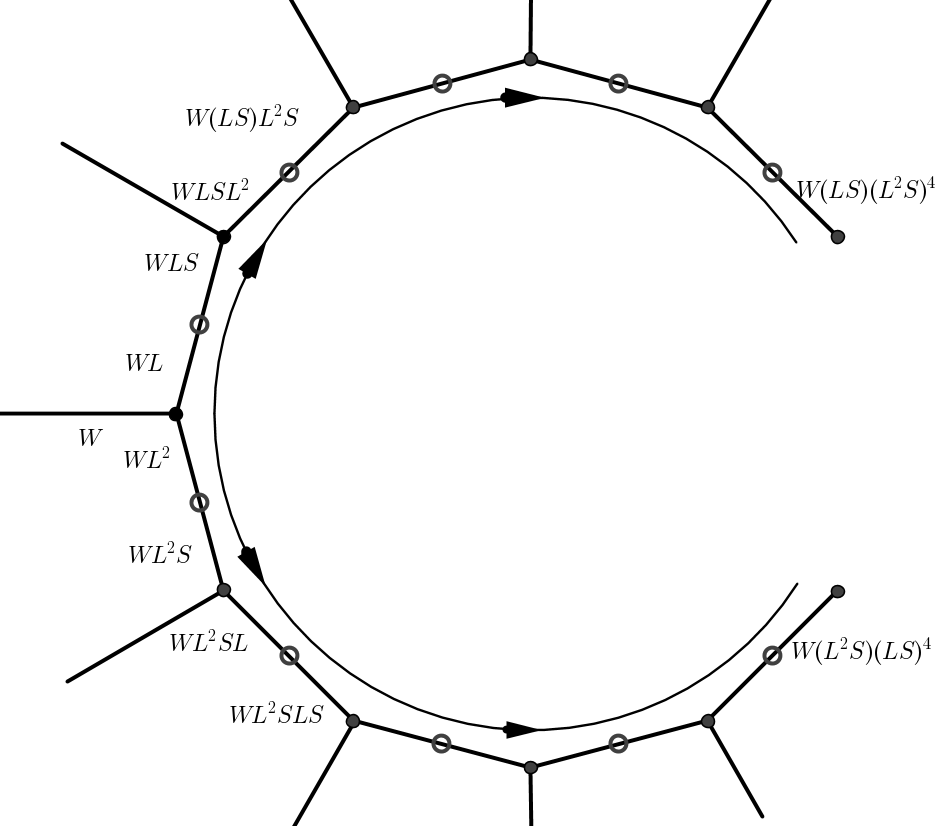}}\\
{\small {\bf Figure.} A pair of rational ends.}\\
\end{center}

On the quotient space $\partial\F_e/\!\sim $ there is the quotient topology induced by the topology on $\partial\F_e$ such that the projection map 
\begin{equation}\label{projectionmap}
 \partial\F_e\longrightarrow \partial\F_e/\!\sim  
 \end{equation}
  is continuous. The quotient space $S^1_e:=\partial\F_e/\!\sim $ is a cyclically ordered topological space under the order relation inherited from $\partial\F_e$.  The equivalence relation is preserved under the canonical homeomorphisms $\partial\F_e \longrightarrow \partial\F_{e'}$ and is also respected by the $\psl$-action. Therefore we have the commutative diagram 
\begin{displaymath}
    \xymatrix{
\partial\F_e  \ar@{->}[r] \ar@{->}[d]  &  \partial\F_{e'}\ar@{->}[d]\\ 
S^1_{e}            \ar@{->}[r]                       &  S^1_{e'}\\
 }
\end{displaymath}
where the horizontal arrows are order-preserving homeomorphisms. 
Moreover, $\psl$ acts by homeomorphisms on $S^1_{e}$, for all $e$.

Now, $\F$ is equipped with a distinguished edge, the one marked $I$, the identity element of $\psl$. So all spaces $S^1_{e}$ are canonically homeomorphic to $S^1_{I}$. 

Any element of $S^1_I$ can be represented by an infinite word in $L$ and $S$. Regrouping occurrences of $LS$ and $L^2S$, any such word $x$ can be written as
\begin{align}\label{reps}
	x&=(LS)^{n_0}(L^2S)^{n_1}(LS)^{n_2}(L^2S)^{n_3}(LS)^{n_4} \cdots \mbox{ or  }\\
	x&=S(LS)^{n_0}(L^2S)^{n_1}(LS)^{n_2}(L^2S)^{n_3}(LS)^{n_4} \cdots,
\end{align}

\noindent where $n_0, n_1 \dots \geq 0$. Since our paths do not have any backtracking we have $n_0\geq0$ and $n_{i}>0$ for $i = 1,2,\cdots$. 
The pairs of words 
\begin{align*}
	(LS)^{n_0}\cdots (LS)^{n_{k}+1} (L^2S)^{\infty} \mbox{ and  } &(LS)^{n_0}\cdots (LS)^{n_{k}}(L^2S) (LS)^{\infty}, & (k\mbox{ even}) \\
	(LS)^{n_0}\cdots (L^2S)^{n_{k}+1} (LS)^{\infty} \mbox{ and  } &(LS)^{n_0}\cdots (L^2S)^{n_{k}}(LS) (L^2S)^{\infty}, & (k\mbox{ odd})  
\end{align*}
correspond to pairs of rational ends and represent the same element of $S_I^1$. For irrational ends this representation is unique. 

The $\psl$-action on $S^1_I$ is then the pre-composition of the infinite word by the word in $L,S$ representing the element of $\psl$. In this picture it is readily seen that this action respects the equivalence relation
$\sim $.

Set $T:=LS$, so that $T(x)=1+x$. Note that
\begin{eqnarray*}
(LS)^n.(L^2S)^m.(LS)^k(x)=(LS)^n.S.[S.(L^2S)^m .S]. S . (LS)^k(x)\\
=(LS)^n.S.[SL^2]^m. S . (LS)^k(x)=(LS)^n.S.[LS]^{-m}. S . (LS)^k(x)\\
=(x+n)\circ (-1/x) \circ (x-m) \circ (-1/x)\circ (x+k) =
n+\cfrac{1}{m+\cfrac{1}{k+x}}
\end{eqnarray*}
Accordingly, define the {\it continued fraction map} $\kappa: S^1_I\rightarrow \hat {\mathbf R} $ by
$$
\kappa(x)=\left\{ \begin{array}{rl}
[n_0,n_1,n_2,\dots]& \mbox{ if } x=(LS)^{n_0}(L^2S)^{n_1}(LS)^{n_2}(L^2S)^{n_3}(LS)^{n_4}\dots\\  
-1/[n_0,n_1,n_2,\dots]& \mbox{ if } x=S(LS)^{n_0}(L^2S)^{n_1}(LS)^{n_2}(L^2S)^{n_3}(LS)^{n_4}\dots
\end{array}\right.
$$\\[-10mm]
\begin{theorem}
	The continued fraction map $\kappa$ is a homeomorphism.
	\label{thm:cfm}
\end{theorem}

\noindent{\it Proof.}
To each rational end $[n_0,n_1,n_2,\dots,n_k,\infty]$ we associate the rational number $[n_0,n_1,n_2,\dots,n_k]$. Likewise for the rational ends in the negative sector. This is an order preserving bijection between the set of equivalent pairs of rational ends and 
$\Q\cup\{\infty\}$.  Now observe that an infinite path is no other than a Dedekind cut and conversely every cut determines a unique infinite path, see  \cite{UZD} for details.\hfill $\Box$

\bigskip
As a consequence, we see that $\kappa$ conjugates the $\psl$-action on $S_I^1$ to its action on  
 $\widehat {\mathbf R}$ by M\"obius maps.
Furthermore, there is a bijection between $\widehat {\Q}$ and the set of pairs of equivalent rational ends. 
Any pair of equivalent rational ends determines a unique {\it rational horocycle}, a bi-infinite left-turning (or right-turning) path.  
Hence, there is a bijection between $\widehat {\Q}$ and the set of rational horocycles. 

\nt{a figure here}

\begin{lemma} \label{correspondences} Given the base edge $I$ of $\F$,\\
(i) There is a natural bijection between the 3-valent vertices of $\F$ and $\widehat {\Q}\setminus\{0,\infty\}$.\\
(ii) There is a bijection between the 2-valent vertices and the Farey intervals $[p/q, r/s]$ with $ps-qr=1$.
\end{lemma}
{\it Proof.}
(i) On each rational horocycle there lies a unique trivalent vertex which is closest to $I$. If we exclude the two horocycles on which $I$ lies, this gives a bijection between the set of trivalent vertices of $\F$ and the set of horocycles. The  map $\kappa$ sends
the two horocycles through $I$ to $0$ and $\infty$. Hence, there is a bijection between the set of trivalent vertices of $\F$ and 
$\widehat {\Q}\setminus \{0, \infty\}$. This bijection sends the trivalent vertex $\{W, WL, WL^2\}$ to $W(1)\in \Q$, where $W$ is a reduced word ending with $S$. (ii) Every 2-valent vertex $\{W, WS\}$ lies exactly on two horocycles, and the set of paths  through $\{W, WS\}$ is sent to the interval $[W(0), WS(0)]$ under $\kappa$. \hfill $\Box$
\sherhh{According to the lemma, there is a correspondence between the one-third of $\psl$ (i.e. the set of cosets of $\langle L\rangle$ = the set of trivalent vertices) and $\widehat {\Q}\setminus \{0, \infty\}$.

On the other hand, note that, given a point at infinity on the unit circle, the torus $\widehat{\R}^1\times \widehat{\R}^1$ can be viewed as the set of intervals on the circle, and it is possible to define a map 
$$
\Jimm: \widehat{\R}^1\times \widehat{\R}^1\to \widehat{\R}^1\times \widehat{\R}^1,
$$
sending intervals to intervals. If $\alpha<\beta$, then 
$$
[\alpha, \beta] \to [\lim_{x\to \alpha^+}\Jimm(x),  \lim_{x\to \beta^-}\Jimm(x) ]
$$
and
$$
[\beta, \alpha] \to [\lim_{x\to \beta^+}\Jimm(x),  \lim_{x\to \alpha^-}\Jimm(x) ]
$$
Note that this map is well-defined at each point (and this holds true if the above trick is applied to a function with only jump discontinuities).
}

\sherh{One may use this lemma to introduce a modular-group structure directly on the set of rationals (i.e. by concatenating paths).}

\paragraph{Periodic paths and the real multiplication set.}
Let $\gamma:=(W_1, W_2, \dots, W_n)$ be a finite path in $\F$. Then the {\it periodization} of $\gamma$ is the path $\gamma^\omega$ defined as
\begin{eqnarray*}
W_1, W_2, \dots, W_n,   W_nW_1^{-1}W_2,     W_nW_1^{-1}W_3,    \dots,    W_nW_{1}^{-1}W_n,  \qquad\qquad\qquad \\ 
\quad (W_nW_{1}^{-1})^2W_2, \dots, (W_nW_{1}^{-1})^2W_n, (W_nW_{1}^{-1})^3W_2, \dots
\end{eqnarray*}
In plain words, $\gamma^\omega$ is the path obtained by concatenating an infinite number of copies of a path representing $W_1^{-1}W_n$, starting at the edge $W_1$. If $W_1^{-1}W_n$ is elliptic, then $\gamma^\omega$ is an infinitely backtracking finite path. If not, 
$\gamma^\omega$ is actually infinite and represents an end of $\F$. We call these {\it periodic ends of $\F$}. 
Thus we have the periodization map
$$
per:\,\,   \{\,finite\, \, hyperbolic \, \, paths\,  \} \longrightarrow \{\, ends\, \, of\, \, \F\,  \}
$$
whose image consists of periodic ends. The map $per$ is not one-to one but its restriction to the set of primitive paths is. 
The modular group action on $\partial\F$ preserves the set of periodic ends and $per$ is $\psl$-equivariant.

Given an edge $e$ of $\F$ and a periodic end $b$ of $\F$, there is a unique path in the class $b$ which starts at $e$. This way the set of periodic ends of $\F$ is identified with the set of eventually periodic paths based at $e$. This set is dense in $\partial\F_e$  and preserved under the canonical homeomorphisms 
between $\partial\F_e$ and $\partial\F_{e'}$. Every periodic end has a unique $\psl$-translate, which is a purely periodic path based at $e$.
Finally, the set of periodic ends descends to a well-defined subset of $S^1_e$. The image of this set under $\kappa$ consists of the set of eventually periodic continued fractions, i.e. the set of real quadratic irrationals. These are precisely the fixed points of the $\psl$-action on $S^1_e$; 
a hyperbolic element $M=(LS)^{n_0}(L^2S)^{n_1}\cdots (LS)^{n_{k}}\in \psl$ fixing the numbers represented by the infinite words
\begin{eqnarray*}
(LS)^{n_0}(L^2S)^{n_1}\cdots (LS)^{n_{k}}(LS)^{n_0}(L^2S)^{n_1}\cdots (LS)^{n_{k}}\dots, \mbox{ and}\\
(LS)^{-n_k}\cdots (L^2S)^{-n_1}(LS)^{-n_{k}}(LS)^{-n_k}\cdots (L^2S)^{-n_1}(LS)^{-n_{k}}\dots. 
\end{eqnarray*}
\sherh{{\bf Exercice.} Show that the fixed points are purely periodic if and only if $M$ is a cyclically reduced word.}
\subsection{The automorphism group of $\F$}
Now let us forget about the ribbon structure of  $\F$.
This gives an abstract graph which we denote  by $|\F|$.
The automorphism group $Aut(|\F|)$  of $|\F|$  contains $Aut(\F)$ as a non-normal subgroup.
It acts by homeomorphisms on $\partial \F$. However, automorphisms of $|\F|$ do not respect the ribbon structure on $\F$ in general and don't induce well-defined self-maps of  $S_I^1$.

 Fix an edge $e$ of $\F$ and denote by $Aut_e(|\F|)$  the group of automorphisms of $|\F|$ that stabilize  $e$. For any pair $e$, $e'$ of edges, $Aut_e(|\F|)$ and $Aut_{e'}(|\F|)$ are conjugate subgroups.
For $n>0$ let $|\F_n|$ be the finite subtree of $|\F|$ containing vertices of distance $\leq n$ from $e$. 
Then $|\F_{n}|\subset |\F_{n+1}|$ forms an injective system with respect to inclusion and $Aut_e(|\F_{n+1}|) \rightarrow Aut_e(|\F_{n}|)$ forms a projective system,
and 
$$
Aut_e(|\F|)=\lim_{\longleftarrow} Aut_e(|\F_{n}|),
$$
showing that $Aut_e(|\F|)$ is a profinite group.  

\paragraph{A description of tree automorphisms as shuffles.} Let us turn back to the ribbon graph $\F$, with the base edge $I\in\psl$. Given any vertex $v$, this base edge permits us to define the subtree (i.e. Farey branch) attached to $\F$ at $v$.

Denote by $V_\bullet$ the set of of degree-3 vertices of $\F$.

The ribbon structure of $\F$ serves as a sort of coordinate system to describe all automorphisms of $|\F|$, as follows.
Given a vertex $v=\{W, WL, WL^2\}$ of type $V_\bullet$ of $\F$, the {\it shuffle} $\sigma_v$ is the automorphism of $|\F|$ which is defined as:
\begin{eqnarray}\label{twist}
\sigma_v: \mbox{edge labeled } M \longrightarrow 
\mbox{edge labeled} \left\{
\begin{array}{ll}
M, & \mbox{ if } M\neq WLX\\ 
WL^2X,& \mbox{ if } M=WLX\\
WLX,& \mbox{ if } M=WL^2X
\end{array}
\right.
\end{eqnarray}
where $M$ and $X$ are reduced words in $S$ and $L$. Thus $\sigma_v$ is the identity away from the Farey branches at $v$, whereas it exchanges the two Farey branches at $v$.
Note that $\sigma_v^2=I$, i.e. the shuffle $\sigma_v$ is involutive.

The automorphism $\sigma_v$ is obtained by permuting the two branches attached at $v$, by shuffling these branches one above the other. 
Beware that $\sigma_v$ is {\it not} the automorphism $\theta_v$ of $|\F|$ obtained by rotating  in the physical 3-space the branches starting at $WL$ and at $WL^2$ around the vertex $v$. We call this latter automorphism a {\it twist}, see below for a precise definition.

We must  stress that the definition of $\sigma_v$ 
requires the ribbon structure of $\F$ as well as a base edge, although $\sigma_v$  is never an automorphism of $\F$. 

Evidently, $\sigma_v$ stabilizes $I$, so one has $\sigma_v\in Aut_I(|\F|)$.

Given an arbitrary (finite or infinite) set $\nu$ of vertices in $V_\bullet$, we inductively define the shuffle 
$\sigma_\nu \in Aut_I(|\F|)$ as follows:
First order the elements of $\nu$ with respect to the distance from the base edge $I$, i.e. set 
$\nu^{(1)}=(v_1^{(1)}, v_2^{(1)},\dots)$, where $v_i^{(1)}\neq v_{i+1}^{(1)}$ and such that 
$
d(v_i^{(1)},I) \leq d(v_{i+1}^{(1)},I)$ for all $i=1,2,\dots.$
Set
$$
\sigma^{(i)}_\nu:=\sigma_{v_i^{(i)}} \mbox{ and } v_j^{(i+1)}=\sigma^{(i)}_\nu(v_j^{(i)}) \mbox{ for } j\geq i+1, \quad 
i=1,2,\dots
$$
Since for any $j>i$, the automorphism $\sigma^{(j)}_\nu$ agrees with  $\sigma^{(i)}_\nu$ on $\F_i$,
the sequence $\sigma^{(i)}_\nu$ converges in $Aut_I(|\F|)$ and we set
$$
\sigma_\nu:=\lim_{i\rightarrow \infty} \sigma^{(i)}_\nu
$$
for the limit automorphism of  $|\F|$. There is some arbitrariness in the initial ordering of $\nu$, concerning its elements of constant distance to the edge $I$, but the limit does not depend on this. The reason is that the shuffles corresponding to those elements commute. 
$\sigma_\emptyset$ is the identity automorphism by definition. 
\begin{theorem}
$Aut_I(|\F|) =\{ \sigma_\nu\, |\, \nu\subseteq V_\bullet(\F)\}$.
\end{theorem}

\noindent{\it Proof.} It suffices to show that the restrictions of shuffles to finite subtrees $|\F_n|$ 
gives the full group $Aut_I(|\F_n|)$. This is easy.
\hfill$\Box$

\paragraph{The automorphism $\sigma_{V_\bullet}$.}
What happens if we shuffle every $\bullet$-vertex of $|\F|$? 
In other words, what is effect of the automorphism $\sigma_{V_\bullet}$ on 
$\partial_I F$?
If $x\in \partial_I F$ is represented by an infinite word in $L$, $L^2$ and $S$, then according to (\ref{twist}), we see that $\sigma_{V_\bullet}$
replaces $L$ by $L^2$ and vice versa. 
In other words, if $x$ is of the form 
$$	
x=(LS)^{n_0}(L^2S)^{n_1}(LS)^{n_2}(L^2S)^{n_3}(LS)^{n_4}\dots, 
$$
then one has
$$
\sigma_{V_\bullet}(x)=
(L^2S)^{n_0}(LS)^{n_1}(L^2S)^{n_2}(LS)^{n_3}(L^2S)^{n_4}\dots
$$
In terms of the contined fraction representations, we get, for $n_0>0$,
$$
\sigma_{V_\bullet}([n_0, n_1,n_2,\dots])=[0, n_1,n_2,\dots]
$$
and for $n_0=0$
$$
\sigma_{V_\bullet}([0, n_1,n_2,\dots])=[n_1,n_2,\dots].
$$
A similar formula also holds when $x$ starts with an $S$, and we get
$$
\sigma_{V_\bullet}(x)=1/x=U(x).
$$ 
Observe (keeping in mind the forthcoming parallelism between $\Jimm$ and $U$) that $U$ is an involution satisfying the equations 
\begin{equation}\label{fesforu}
U=U^{-1}, \quad US=SU, \quad
LU=UL^2.
\end{equation}
Here, $U$, $L$ and $S$ are viewed as operators acting on the boundary $\partial_I\F$.
These equations can be re-written in the form of  functional equations as below:
$$
U(Ux)=x, \quad U(-\frac{1}{x})=-\frac{1}{Ux}, \quad U(\frac{1}{1-x})=1-\frac{1}{Ux}.
$$
One may derive other functional equations from these, i.e. $UT=ULS=L^2US=L^2SU$ is written as
$$
U(1+x)=\frac{Ux}{1+Ux}.
$$
One may consider $U$ as an ``{orientation-reversing}" automorphism of the ribbon tree $\F$. 
In the same vein, $U$ is a homeomorphism of 
$\partial_I\F$ which reverses its canonical ordering. It is the sole
element of $Aut_I(|\F|)$ which respects the equivalence $\sim $
(unlike $\Jimm$) and hence $U$ acts by homeomorphism on $S^1_I$.
In fact, taken as a map on $E(\F)=\psl$, the involution $U$ is an automorphism of $\psl$, the one defined by 
$$
U:\left(\begin{matrix} S\\L\end{matrix}\right) \to 
\left(\begin{matrix} S\\L^2\end{matrix} \right)
$$
The automorphism group of $\psl$ is generated by the inner automorphisms and $U$: this is the group $\pgl$.

\paragraph{Twists.} Let $v\in V_\bullet$ and let 
$\nu_v$ be the set of all vertices (including $v$) on the Farey branch (with reference to the edge $I$) of $\F$ at $v$. Then the {\it twist} of $v$
 is the automorphism of $\F$ defined by 
$$
\theta_v:=\sigma_{\nu_v}.
$$
In words, $\theta_v$ is the shuffle of every vertex of the Farey branch at $v$.
As in the case of shuffles, for any subset $\mu\subset V_\bullet$, 
one may define the twist $\theta_\mu$ as a limit of 
convergent sequence of individual twists.

\begin{lemma}\label{twisttoshuffle}
Let $v$ be a vertex in $V_\bullet$ and $v'$, $v''$ its two children 
(with respect to the ancestor $I$). Put $\mu:=\{v,v',v''\}$.  
Then
$
\sigma_v=\theta_\mu.
$ 
\end{lemma}
Hence, shuffles can be expressed in terms of twists and vice versa. This proves:
\begin{theorem}
$Aut_I(|\F|) =\{ \theta_\mu\, |\, \mu\subseteq V_\bullet(\F)\}$.
\end{theorem}

\begin{notatdef}\label{notatdef}
Let $v^*$ be the trivalent vertex incident to the base edge $I$, i.e. $v^*:=\{I,L,L^2\}$ and set 
$V_\bullet^*:=V_\bullet\backslash \{v^*\}$. 
We denote the special automorphism $\theta_{V_\bullet^*}$ 
by  $\Jimm_\F$.  In words, this 
is the automorphism of $|\F|$ obtained by 
twisting every trivalent vertex except $v^*$. 
\end{notatdef}
This $\Jimm_\F$ is the same automorphism obtained by 
shuffling every other trivalent vertex, such that  $v^*$ is not shuffled. More precisely, the following consequence of Lemma \ref{twisttoshuffle} holds:
\begin{lemma}
One has 
$
\theta_{V_\bullet^*}=\Jimm_\F=\sigma_J,
$
where $J\subset V_\bullet$ is the set of degree-3 vertices, whose distance to the base edge is an odd number. 
\end{lemma}
Obviously, $\Jimm_\F$ induce an involutive homeomorphism of $\partial_I\F$ 
but violates its ordering.  We denote this homeomorphism by $\Jimm_{\partial\F}$. We must emphasize that 
$\Jimm_{\F}$ and $\Jimm_{\partial\F}$  are perfectly well-defined mappings on their domain of definition.
They don't exhibit such things as the two-valued behaviour of $\Jimm_{\R}$ at rationals. This two-valued behavior is a consequence of the fact that $\Jimm_{\partial\F}$ do not respect the relation $\sim $. 

To see the effect of $\Jimm_{\partial\F}$ on $x\in\partial_I\F$, assume 
$$	
x=S^\epsilon (LS)^{n_0}(L^2S)^{n_1}(LS)^{n_2}(L^2S)^{n_3}(LS)^{n_4}\dots, \, (n_0\geq 0, \, n_i>0 \mbox{ if } i>0, \,\epsilon\in\{0, 1\}).
$$
We may represent $x$ by a string of 0's and 1's (0 for $L$ and 1 for $L^2$):
\begin{eqnarray*}
x=S^\epsilon
\underbrace{00\dots0}_{n_0}
\underbrace{11\dots1}_{n_1}
\underbrace{00\dots0}_{n_2}
\underbrace{11\dots1}_{n_3}
\dots, \, 
\end{eqnarray*}
Then rationals are represented by the eventually constant strings where for any finite string $a$, the strings $a0111\dots$ and $a1000\dots$ represent the same rational\footnote{The two representations of the number $0\in \R$ are
$1^\omega$ and $-0^\omega$ and the  two representations of $\infty\in \widehat\R$ are $0^\omega$ and $-1^\omega$.}.

Let $\phi:=(01)^\omega$ be the zig-zag path to infinity, and set $\phi^*:=\lnot \phi=(10)^\omega$.
Then 
$$
\Jimm_{\partial\F}(x)= \begin{cases} 
a \xor \phi,& x=a\in \{0,1\}^\omega,\\
S(a \xor \phi*),& x=Sa, \, a\in \{0,1\}^\omega.
\end{cases}
$$
where $\xor$ is the operation of term-wise exclusive or (XOR). The involutivity of $\Jimm_{\partial\F}$ then stems from the reversibility of the disjunctive or: $(p\xor q) \xor q=p$.

\sherh{xoring (or performing term-wise logical operations) with a fixed infinite string defines boundary homeomorphisms 
induced by special kind of tree automorphisms, which can be called radially symmetric.}

\sherh{One may take a fixed bi-infinite string with a section, $\alpha^*|\alpha$ and xor with this 
$$
\Jimm_{\partial\F}= \begin{cases} 
a \xor \alpha,& x=a\in \{0,1\}^\omega,\\
S(a \xor \alpha*),& x=Sa, \, a\in \{0,1\}^\omega.
\end{cases}
$$
Example $\alpha^*|\alpha=(100)^*|(001)^*$. Something like this.
These boundary homeomorphisms comes from the tree automorphisms, those which shuffle the 1's on the string $\alpha^*|\alpha$. These are precisely the ``radially symmetric" automorphisms (-twists or -shuffles) of the tree.
}
\medskip\noindent
{\bf Some Examples.}
The string $0(0011)^\omega$ corresponds to the continued fraction 
$[1,2,2,2,\dots]$, which equals $\sqrt{2}$. One has
$$
\begin{array}{l|cccccccccccc|c}
\partial\F&&\multicolumn{10}{c}{\{0,1\}^\omega}&&\widehat{\R}\\
\hline
x 
&0&1&1&0&0&1&1&0&0&1&1&0\dots &\sqrt{2}\\
\phi
&0&1&0&1&0&1&0&1&0&1&0&1\dots&\Phi\\
\Jimm(x)
&0&0&1&1&0&0&1&1&0&0&1&1\dots&1+\sqrt{2}
\end{array}
$$
As for the value of $\Jimm_\R$ at $\infty$, one has
$$
\begin{array}{l|cccccccccccc|c}
\partial\F&&\multicolumn{10}{c}{\{0,1\}^\omega}&&\widehat{\R}\\
\hline
x 
&0&0&0&0&0&0&0&0&0&0&0&0\dots &\infty\\
\phi
&0&1&0&1&0&1&0&1&0&1&0&1\dots&\Phi\\
\Jimm(x)
&0&1&0&1&0&1&0&1&0&1&0&1\dots&\Phi
\end{array}
$$
The two values $\Jimm_\R$ assumes at the point $1$ are found as follows:
$$
\begin{array}{l|cccccccccccc|c}
\partial\F&&\multicolumn{10}{c}{\{0,1\}^\omega}&&\widehat{\R}\\
\hline
x 
&0&1&1&1&1&1&1&1&1&1&1&1\dots &1\\
\phi
&0&1&0&1&0&1&0&1&0&1&0&1\dots&\Phi\\
\Jimm(x)
&0&0&1&0&1&0&1&0&1&0&1&0\dots&1+\Phi
\end{array}
$$
$$
\begin{array}{l|cccccccccccc|c}
\partial\F&&\multicolumn{10}{c}{\{0,1\}^\omega}&&\widehat{\R}\\
\hline
x 
&1&0&0&0&0&0&0&0&0&0&0&0\dots &1\\
\phi
&0&1&0&1&0&1&0&1&0&1&0&1\dots&\Phi\\
\Jimm(x)
&1&1&0&1&0&1&0&1&0&1&0&1\dots&1/(1+\Phi)
\end{array}
$$
Conversely, one has $\Jimm_\R(1+\Phi)=\Jimm_\R(1/(1+\Phi))=1$, 
illustrating the two-to-oneness of $\Jimm_\R$ on the set of noble numbers. 

\bigskip
The noble numbers  are the $\psl$-translates of the golden section $\Phi$.
They correspond to the eventually zig-zag paths in $\partial \F$,
and represented by strings terminating with $(01)^\omega$. 
From the $\xor$-description, it is clear that $\Jimm_{\partial\F}$ sends those strings to rational (i.e. eventually constant) strings and vice versa.

Since the equivalence $\sim $ is not respected by $\Jimm_{\partial\F}$,
it does not induce a homeomorphism of  $\widehat{\R}$, not even a well-defined map. Nevertheless, 
if we ignore the pairs of rational ends, then the remaining equivalence classes are singletons and $\Jimm_{\partial\F}$ restricts to a well-defined map 
$$
\Jimm_\R:\widehat{\R}\setminus\widehat{\Q}
\to \widehat{\R}
$$
If we also ignore the set of ``{noble paths}", i.e. the set
$\pgl \widehat{\Q}=\Jimm^{-1}(\Q)$, then we obtain an involutive bijection
$$
\Jimm_\R: \widehat{\R}\setminus(\widehat{\Q}\cup \pgl \widehat{\Q})\to
\widehat{\R}\setminus(\widehat{\Q}\cup \pgl \widehat{\Q})
$$
\paragraph{The functional equations.}
Since twists and shuffles do not change the distance to the base, 
we have our first functional equation:
\begin{lemma}
The $|\F|$-automorphisms $\Jimm_\F$  and $U=\sigma_{V_\bullet}$ commute, i.e. $\Jimm_\F U=U\Jimm_\F$. Hence,
$\Jimm_{\partial\F}U=U\Jimm_{\partial\F}$ and 
whenever $\Jimm_\R$ is defined, one has 
\begin{equation}\label{jimmu}
\Jimm_\R U= U\Jimm_\R \iff \Jimm_\R\bigl(\frac{1}{x}\bigr)=\frac{1}{\Jimm_\R(x)}
\end{equation}
\end{lemma}
In fact, $\Jimm_\F U$ is the automorphism of $\F$ which shuffles every other vertex, starting with the vertex $v^*$. 
In other words, it shuffles those vertices which are not shuffled by $\Jimm_\F$. We denote this automorphism by 
$\Jimm_\F^*$.
Note that, in terms of the strings, the operation $U$ is nothing but the term-wise negation:
$$
Ua=\lnot a,
$$ 
and the lemma merely states the fact that $\lnot (\phi \xor a)=\phi \xor \lnot a$. The boundary homeomorphism induced by the 
automorphism $\Jimm_\F^*$ is thus the map $\Jimm_{\partial \F}^*$ which xors with the string $\phi^*:=(10)^\omega$.

Now, consider the operation $S$. If $x=a\in \{0,1\}^\omega$, then one has $\Jimm_{\partial \F} (Sx)=$
\begin{eqnarray*}
 = S(a \xor \phi^*) 
=S(a\xor \lnot \phi) =S(\lnot(a\xor\phi))=S(U(a\xor\phi)) = SU\Jimm_{\partial \F} (x)=V\Jimm_{\partial \F} (x)
\end{eqnarray*}
The same equality holds if $x=Sa$, and we get our second functional equation:
\begin{lemma}
The $\partial_I\F$-homeomorphisms $\Jimm_{\partial\F}$  and $S$ satisfy $\Jimm_{\partial\F}S=V\Jimm_{\partial\F}$. 
Hence, whenever $\Jimm_\R$ is defined, one has 
\begin{equation}\label{jimmux}
\Jimm_\R S= V\Jimm_\R \iff \Jimm_\R\bigl(-\frac{1}{x}\bigr)=-{\Jimm_\R(x)}
\end{equation}
\end{lemma}

Now we consider the operator $L$. Suppose that $x=Sa$.
Then $Lx= LSa=0a$. Hence, noting that $\phi=(01)^\omega=0 (10)^\omega=0\phi^*$ we have
\begin{eqnarray*}
 \Jimm_{\partial \F} (Lx) = 0a \xor \phi=0a\xor 0\phi^*=0(a\xor \phi^*)= LS(a\xor \phi^*) 
 = L\Jimm_{\partial \F} (x).
\end{eqnarray*}
Another possibility is that $x=0a$.
In other words, $x$ starts with an $L$. 
Then $Lx$ starts with an $L^2$, i.e. $Lx= 1a$. Hence, 
\begin{eqnarray*}
 \Jimm_{\partial \F} (Lx) = 1a \xor \phi=1a\xor 0\phi^*=1(a\xor \phi^*),\\
  \Jimm_{\partial \F}(x)=0a \xor \phi=0a\xor 0\phi^*=0(a\xor \phi^*) 
\implies L \Jimm_{\partial \F}(x)= 1(a\xor \phi^*).
\end{eqnarray*}
Finally, if $x=1a$, then
$Lx$ starts with an $S$, i.e. $Lx= Sa$. Hence, 
\begin{eqnarray*}
 \Jimm_{\partial \F} (Lx) = S(a \xor \phi^*),\\
  \Jimm_{\partial \F}(x)=1a \xor \phi=1a\xor 0\phi^*=1(a\xor \phi^*) 
\implies L \Jimm_{\partial \F}(x)= S(a\xor \phi^*).
\end{eqnarray*}
Whence the third functional equation:
\begin{lemma}
The $\partial_I\F$-homeomorphisms $\Jimm_{\partial\F}$  and $L$ satisfy $\Jimm_{\partial\F}L=L\Jimm_{\partial\F}$. 
Hence, whenever $\Jimm_\R$ is defined, one has 
\begin{equation}\label{jimmux}
\Jimm_\R L= L\Jimm_\R \iff \Jimm_\R\bigl(1-\frac{1}{x}\bigr)=1-\frac{1}{\Jimm_\R(x)}
\end{equation}
\end{lemma}
Since $U$, $S$ and $L$ generate the group $\pgl$, these three functional equations  forms a complete set, from which the rest can be deduced. For example, 
$$
T=LS\implies \Jimm_{\partial \F}T=LV\Jimm_{\partial \F}\iff
\Jimm_\R(1+x)=1+\frac{1}{\Jimm_\R(x)}
$$
This also shows that $\Jimm_\R$ acts as the desired outer automorphism of $\pgl$. 
\paragraph{Other forms of the functional equations.}
It is possible to derive many alternative forms of the functional equations.
Here we record some of them. We leave the task of verifying to the reader. 
(As usual we drop the subscript $\R$ from $\Jimm_\R$ to increase the readability):
$$
\Jimm\left(1-\frac{1}{x}\right)=1-\frac{1}{\Jimm(x)}, \quad 
\Jimm\left(\frac{x}{x+1}\right)+\Jimm\left(\frac{1}{x+1}\right)=1.
$$
\smallskip\begin{center}\label{twovariable}
{\fbox{\begin{minipage}{10cm}
$$
\Jimm(x)=y  \iff \Jimm(y)=x
$$

$$
xy=1\iff \Jimm(x)\Jimm(y)=1
$$

$$
x+y=0\iff \Jimm(x)\Jimm(y)=-1
$$

$$
x+y=1\iff \Jimm(x)+\Jimm(y)=1
$$

$$
\frac{1}{x}+\frac{1}{y}=1\iff \frac{1}{\Jimm(x)}+\frac{1}{\Jimm(y)}=1
$$

\vspace{1mm}
\end{minipage}}}
\end{center}

\section{Diverse facts about $\Jimm$}
\subsection{$\Jimm$ as a limit of piecewise-$\pgl$ functions.}
In $\F$, denote by  
$V_\bullet(n)$ the set of vertices of 
distance $\leq n$ to the base (excluding $v^*$), and let $V_\bullet^\prime(n)$ be the set of vertices of odd distance $\leq n$ to the base. Then
$$
\Jimm_\F=
\lim_{n\to\infty} \theta_{V_\bullet(n)}
=\lim_{n\to\infty} \sigma_{V_\bullet^\prime(n)}.
$$
The maps $\theta_{V_\bullet(n)}$ and $\sigma_{V_\bullet^\prime(n)}$ 
induce a sort of finitary 
projective interval exchange maps on $\widehat{\R}$
and thus $\Jimm_\R$ can be written as a limit of such functions. 
Below we draw the twists $\theta_{V_\bullet(n)}$ for $n=2\dots 5$. 
\sherh{These are obtained by xoring with the sequences $(01)^n1^\omega$ or $(01)^n0^\omega$. The non-defined values 0 and $\infty$ can be thought of those rationals having nothing to be xored with.}
\begin{center}\label{jimmiterate}
\noindent
{\includegraphics[scale=.2]{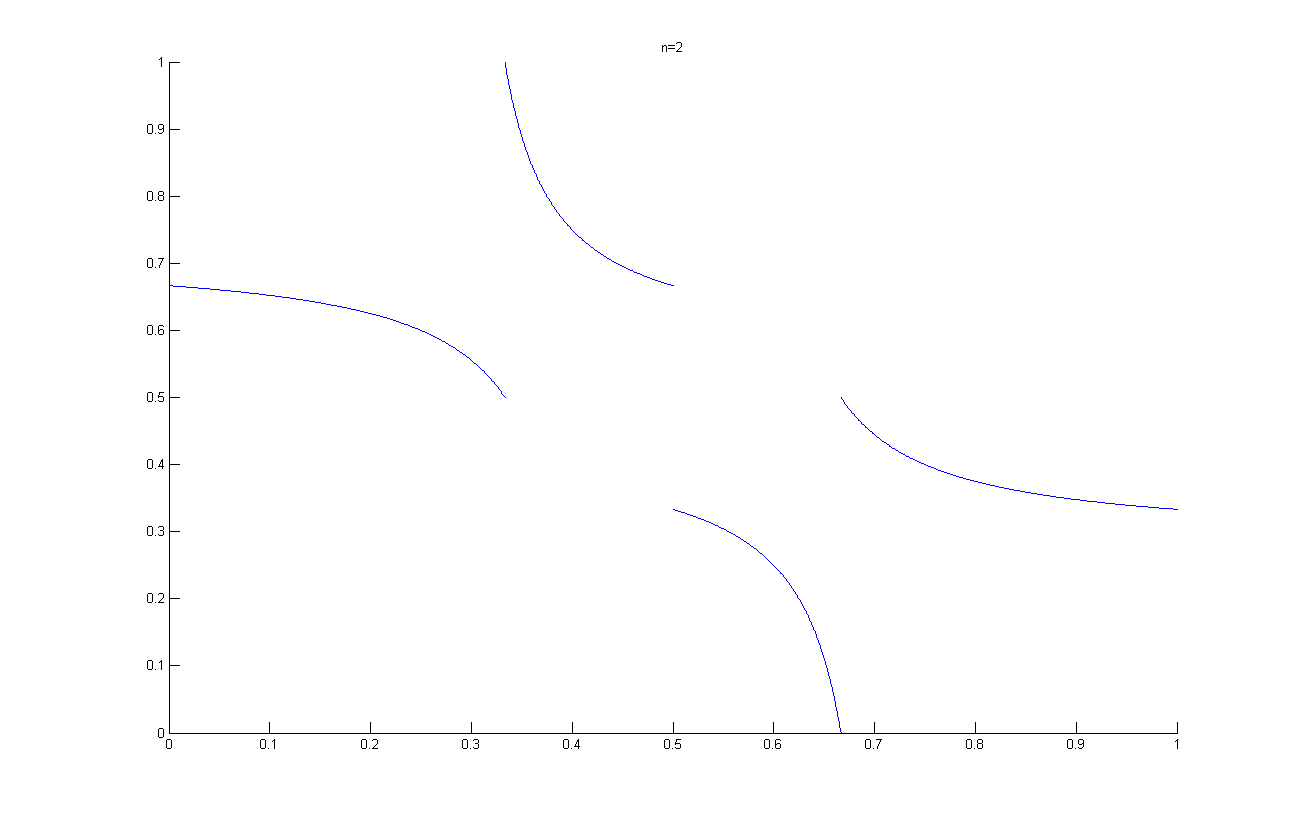}}
{\includegraphics[scale=.2]{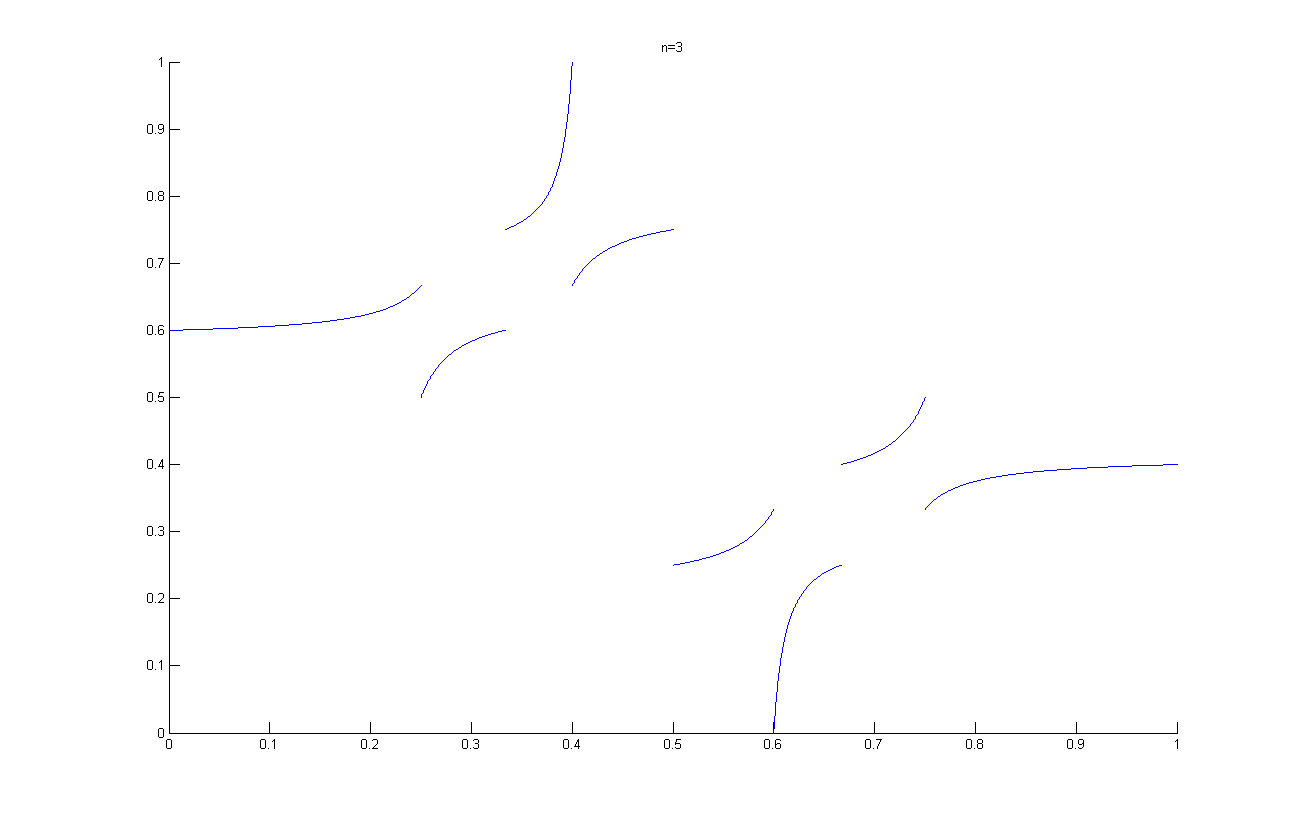}}
{\includegraphics[scale=.2]{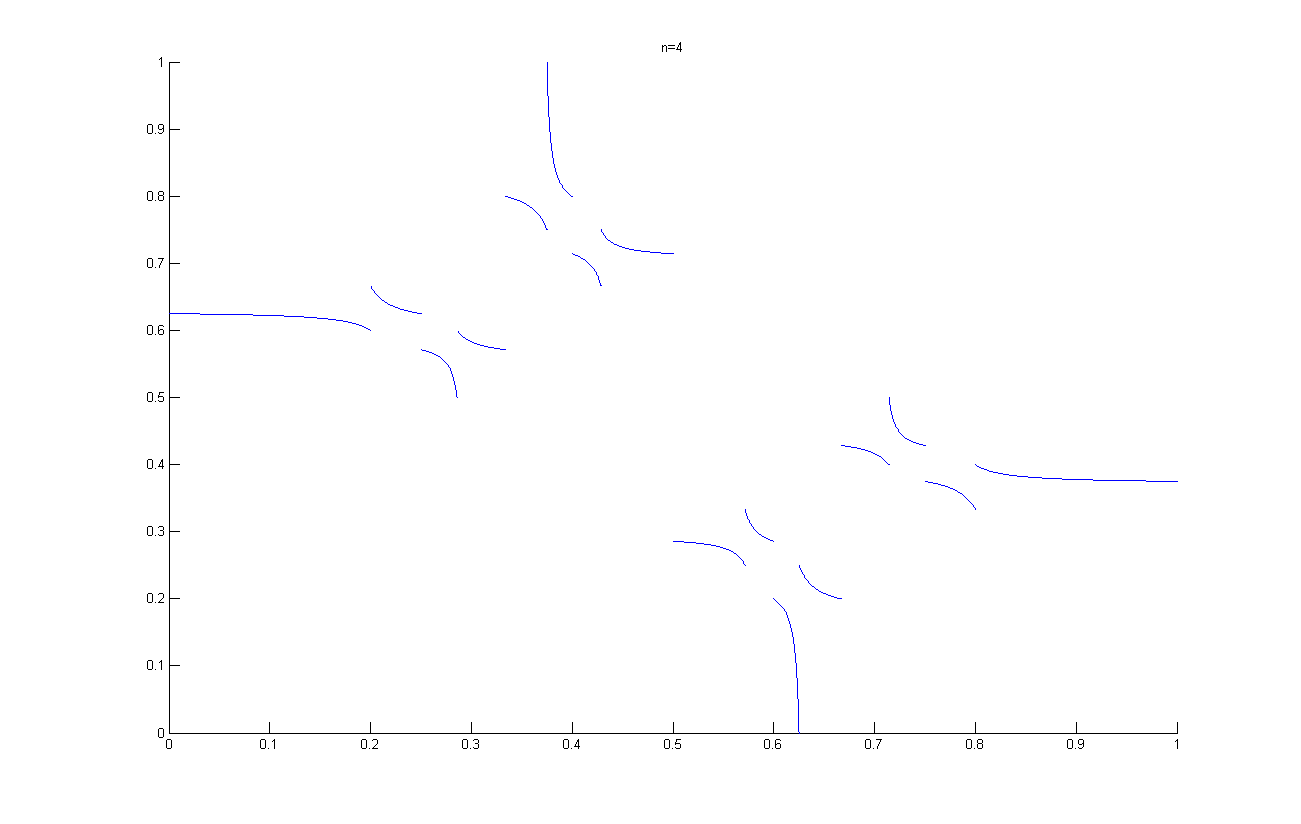}}
{\includegraphics[scale=.2]{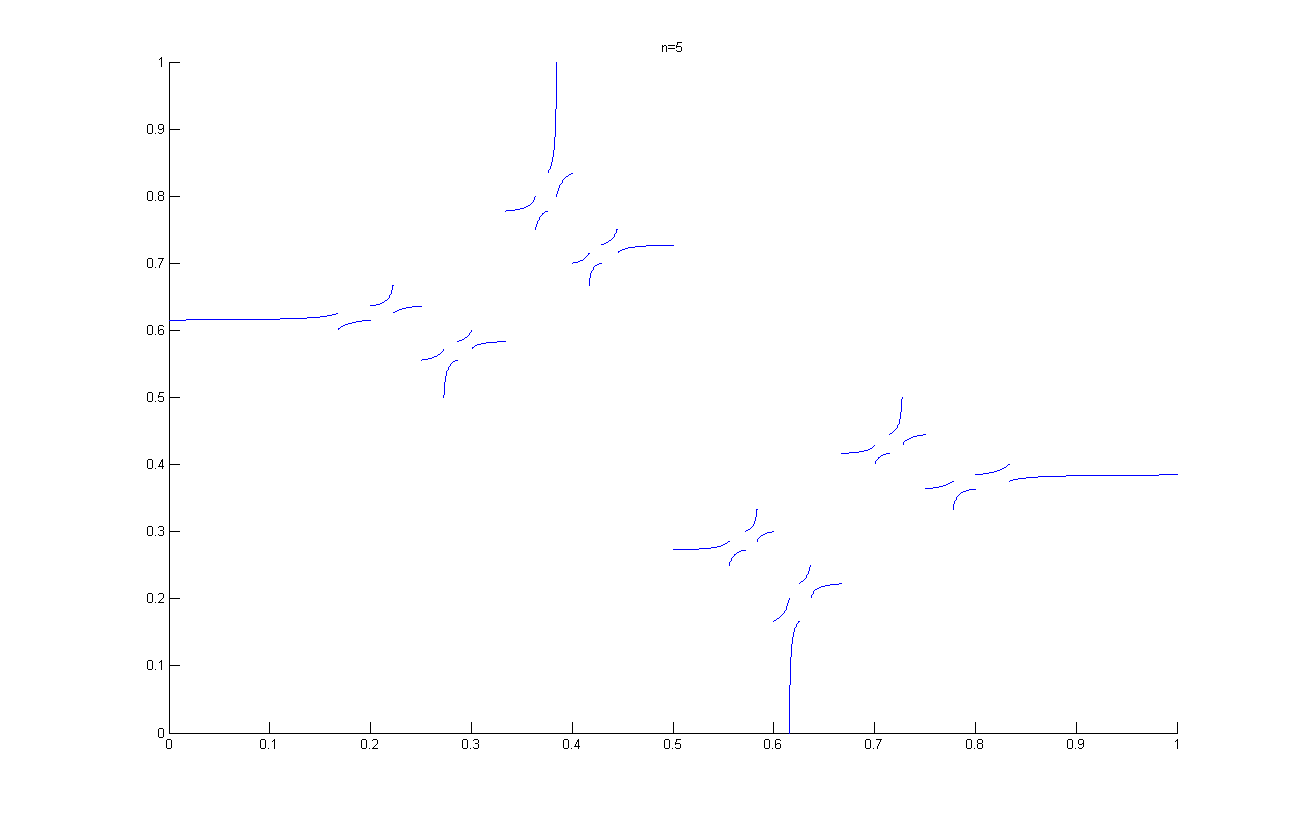}}\\
{\small {\bf Figure.} The plot of $\theta_{\F_n^*}$ on the interval $[0,1]$, for $n=2,3,4,5$.}\\
\end{center}

\subsection{$\Jimm$ on rationals}\label{jimmq}
In virtue of Lemma \ref{correspondences}, there is a 1-1 correspondence between the set of trivalent vertices and the set 
$\Q\setminus\{0,1\}$. Since any element of $Aut_I(\F)$ defines a bijection of $V_\bullet(\F)$, we see that every automorphism of $\F$ that fixes $I$, defines a unique bijection of $\Q\setminus\{0,1\}$. In particular, this is the case with $\Jimm_\F$. We extend this involution to $\Q\cup\{\infty\}$ by sending $0\leftrightarrow \infty$ and we denote the resulting
bijection with $\Jimm_\Q$. This involution satisfies all the functional equations satisfied by $\Jimm_\R$. One has $\Jimm_\Q(1)=1$, and for $x>0$, its values can be computed by using the functional equations 
$\Jimm_\Q(1+x)=1+1/\Jimm_\Q(x)$ and $\Jimm_\Q(1/x)=1/\Jimm_\Q(x)$. It tends to $\Jimm_\R$ at irrational points and in fact $\Jimm_\R$ can be defined as this limit. However, it must be emphasized that $\Jimm_\Q(x)$ is not the restriction of $\Jimm_\R$ to $\Q$; this latter function is by definition two-valued at rationals\footnote{It might have been convenient to declare the values of $\Jimm_\R$ at rational arguments to be given by $\Jimm_\Q(x)$, so that $\Jimm_\R$ would be a well-defined function everywhere (save 0 and $\infty$). However, we have chosen to not to follow this idea, for the sake of uniformity in definitions.}.

We show  in \cite{lebesguesym} that  $\Jimm_\Q$ ``commutes" with Lebesgue measure in a certain sense.

The involution $\Jimm_\Q$ conjugates the multiplication (denoted by $\odot$) on $\Q$ to an operation with 1 as its identity, and such that the inverse of $q$ is $1/q$. The addition (denoted by $\oplus$) is conjugated to an operation with $\infty$ as its neutral element and such that the additive inverse of $q$
is $-1/q$, i.e. $\ominus q=-1/q$.

\subsection{Fibonacci sequence}
The Fibonacci sequence is defined by the recurrence 
$
F_0=0, \quad F_1=1,\mbox{ and } F_{n}=F_{n-1}+F_{n-2} \mbox{ for } n\in \Z.
$
One has then $F_{-n}=(-1)^{n+1}F_n$ and
$$
\widetilde{T}(x)=1+\frac{1}{x}\implies
\widetilde{T}^n=\frac{F_{n+1}x +F_n}{F_nx+F_{n-1}} \quad (n\in \Z).
$$
The following lemma is an easy consequence of the functional equations.
\begin{lemma}\label{easyconsequence}
Let $x$ be an irrational number. Then\\
(i) 
$$
\Jimm(1+x)=1+\frac{1}{\Jimm(x)} \iff \Jimm(Tx)=\widetilde{T}\Jimm(x).
$$
(ii) 
$$
\Jimm(n+x)=\widetilde{T}^n\Jimm(x)=\frac{F_{n+1}\Jimm(x) +F_n}{F_n \Jimm(x)+F_{n-1}} \quad (n\in \Z).
$$
(iii) 
$$
\Jimm([1,x])=\Jimm(1+\frac{1}{x})=1+\Jimm(x) \implies \Jimm([1_n,x])=n+x \quad (n\in \Z).
$$
(iv)
$$
\Jimm([n,x])=\Jimm(n+1/x)=\frac{F_{n+1}+F_n\Jimm(x) }{F_n +F_{n-1}\Jimm(x)} \quad (n\in \Z).
$$
\end{lemma}
 \subsection{Continuity}
Since the involution $\Jimm_{\partial\F}: \partial\F\to \partial\F$
is a homemorphism and since the subspaces 
$$
\R\setminus \Q \simeq \partial\F\setminus \{\mbox{rational paths}\}
$$ 
are also homeomorphic, the function  $\Jimm_\R$ is continuous at irrational points:
\begin{theorem}\label{jimmcontini}
The function $\Jimm_\R$ on $\R\setminus \Q$ is continuous. 
\end{theorem}
Moreover, recall that there is a canonical ordering on $\partial\F$
inducing on $\widehat{\R}$ its canonical ordering compatible with its topology.
So the notion of lower and upper limits exists on $\partial\F$
and coincides with lower and upper limits on $\widehat{\R}$.
This shows that the two values that $\Jimm_\R$ assumes on rational arguments are nothing but the limits
$$
\Jimm(q)^-:=\lim_{x\to q^-} \Jimm(x) \mbox{ and }
\Jimm(q)^+:=\lim_{x\to q^+} \Jimm(x) 
$$
By choosing one of these values coherently, one can make $\Jimm$ an everywhere upper (or lower) continuous function. Note however that it will (partially) cease to satisfy the functional equations at rational arguments. One has, for irrational $r$ and rational $q$,
$$
\lim_{q\rightarrow r} \Jimm_\Q(q)=\lim_{q\rightarrow r} \Jimm_\R(q)= \Jimm_\R(r),
$$
no matter how we choose the values of $\Jimm_\R(q)$. In fact $\Jimm$ is almost everywhere differentiable with derivative vanishing almost everywhere, see  \cite{jimmarxiv}.

\subsection{Action on the quadratic irrationals}
Since $\Jimm_\F$ sends eventually periodic paths to eventually periodic paths, 
$\Jimm_\R$ preserves the real-multiplication set. This is the content of our next result:

\begin{theorem}\label{orbits}
$\Jimm$ defines an involution of the set of real quadratic irrationals, and it respects the $\pgl$-orbits.
\end{theorem}

Beyond this theorem, we failed to detect any further arithmetic structure or formula relating the quadratic number to its $\Jimm$-transform.  
All we can do is to give some sporadic examples, which amounts to exhibit some special real quadratic numbers with known continued fraction expansions. We shall do this below. Before that, however, note that the equation below is solvable for every $n\in \Z$:
$$
\Jimm x=n+ x\implies  x=\Jimm\Jimm x=\widetilde{T}^n\Jimm x\implies  x=\widetilde{T}^n( x+n)
$$
$$
\frac{F_{n+1}( x+n)+F_n }{F_n( x+n) +F_{n-1}}= x=\frac{F_{n+1} x+F_n+nF_{n+1} }{F_n x^n +F_{n-1}+nF_n}
$$
$$
\implies F_n x^2+(F_{n-1}+nF_n-F_{n+1}) x -(F_n+nF_{n+1})
$$
More generally, the equation $\Jimm x=M x$ is solvable for $M\in\pgl$. Indeed one has
$$
\Jimm x=M x\implies 
 x=\Jimm\Jimm x=(\Jimm M)(\Jimm x)=
(\Jimm M) M x,
$$
and the solutions are the fixed points of $(\Jimm M) M$. 
These points are precisely the words represented by the infinite path
\begin{eqnarray}\label{fixed}
x=(\Jimm M) M(\Jimm M) M(\Jimm M) M(\Jimm M) M\dots
\end{eqnarray}
where we assume that both $M$ and $\Jimm M$ are expressed as words in $U$ and $T$.
Note that these words are precisely the points whose conjugacy classes remain stable under $\Jimm$.
$$
\Jimm x=M(\Jimm M) M(\Jimm M) M(\Jimm M) M(\Jimm M) M\dots
$$
Hence the $\pgl$-orbits of the points in (\ref{fixed}) remain stable under $\Jimm$. Since the set of these orbits is the moduli space of real lattices, we deduce:
\begin{proposition}
The fixed points of the involution $\Jimm$ on the moduli space of real lattices are precisely the points (\ref{fixed}) .
\end{proposition}
For example, if $\Jimm M=\widetilde{T}^{k+1}$  then  $M=T^{k+1}$, and the point (\ref{fixed})  is nothing but the point 
$[\overline{1_k,k+2}]$ mentioned in the introduction. 

\subsection{Examples}
\paragraph{Example.} In  \cite{manfred} it is shown that 
$$
x=[0;{\overline{1_{n-1},a}}]=\frac{a}{2}\left(\sqrt{1+4\frac{aF_{n-1}+F_{n-2}}{a^2F_n}}-1\right)
$$
Hence
$
\Jimm(x) = [0;n,\overline{1_{a-2},n+1}], 
$
and  so $\Jimm(x) = 1/(n+y)$, where
$$
y=[0;\overline{1_{a-2},n+1}]=\frac{n+1}{2}\left(\sqrt{1+4\frac{(n+1)F_{a-2}+F_{a-3}}{(n+1)^2F_{a-1}}}-1\right)
$$

\paragraph{Example.} (From Einsiedler \& Ward \cite{ward}, Pg.90, ex. 3.1.1)
This result of McMullen from~\cite{mcmullen} illustrates how $\Jimm$ behaves on one real quadratic number field.
$\Q(\sqrt{5})$ contains infinitely many elements with a uniform bound on their partial quotients, since 
$
[\overline{1_{k+1},4,5,1_k,3}]\in \Q(\sqrt{5}), \quad \forall k=1,2,\dots.
$
Routine calculations shows that the transforms 
$
\Jimm([\overline{1_{k+1},4,5,1_k,3}])=[k+2, \overline{2,2,3,k+2,1,k+3}]
$
lies in different quadratic number fields for different $k$'s. Hence, not only $\Jimm$ does not preserve the property of ``{belonging to a certain quadratic number field}", it sends elements from one quadratic number field to different quadratic number fields.
The following result provides even more examples of this nature: 
\begin{theorem} (McMullen~\cite{mcmullen})
For any $s>0$, the periodic continued fractions
$$
x_m=[\overline{(1,s)^m,1,s+1,s-1, (1,s)^m,1,s+1,s+3)}]
$$
lie in $\Q(\sqrt{s^2+4s})$ for any $m\geq 0$.
\end{theorem}
\begin{theorem}(Fishman and Miller \cite{fishman})
One has the following
$$
\Phi^k=
\begin{cases}
[\overline{F_{k+1}+F_{k-1}}],& \mbox{\it if k is even},\\
[F_{k+1}+F_{k-1}-1, \overline{1,F_{k+1}+F_{k-1}-2}]& \mbox{\it if k is odd}.
\end{cases}
$$
\end{theorem}
\begin{corollary}
The continued fraction of the kth power of the golden section is
$$
\Jimm(\Phi^k)=
\begin{cases}
[1_{F_{k+1}+F_{k-1}-1},\overline{2, 1_{F_{k+1}+F_{k-1}-2}}],& \mbox{\it if k is even},\\
[1_{F_{k+1}+F_{k-1}-2},\overline{3,1_{F_{k+1}+F_{k-1}-4}}]& \mbox{\it if k is odd}.
\end{cases}
$$
\end{corollary}

\paragraph{Example.}  Suppose $\alpha=p+\sqrt{q}$ be a quadratic irrational with $p,q\in \Q$. Then $\alpha^*=1/\alpha$ provided
$$
|\alpha|^2=p^2-q=1 \implies q=p^2-1.
$$
Suppose $\alpha$ is of this form, i.e. $\alpha=p+\sqrt{p^2-1}$ and suppose $\Jimm(\alpha)=x+\sqrt{y}$. 
Then since 
$$
x-\sqrt{y}=\Jimm(\alpha)^*=\Jimm(\alpha^*)=\Jimm(1/\alpha)=1/\Jimm(\alpha)=1/(x+\sqrt{y})
\implies x^2-y=1,
$$
and we conclude that $\Jimm(\alpha)$ is again of the form $x+\sqrt{x^2-1}$. 

On the other hand, suppose $\Jimm(\sqrt{q})=x+\sqrt{y}$. Then 
$$
\Jimm(\sqrt{q})^*=\Jimm(\sqrt{q}^*)=\Jimm(-\sqrt{q})=-1/\Jimm(\sqrt{q})\implies
|\Jimm(\sqrt{q}^*)|=x^2-y=-1
$$
Hence, $\Jimm(\sqrt{q})$ is of the form $x+\sqrt{1+x^2}$. Conversely, 
if $\alpha$ is of the form $p+\sqrt{p^2+1}$, then $\Jimm(\alpha)$ is a quadratic surd.

\sherh{because
$$
a-2+\frac{1}{n+1+\frac{1}{y}}=y\implies a-2+\frac{y}{(n+1)y+1}=y \implies 
$$$$
(a-2)[(n+1)y+1]=(n+1)y^2\implies (n+1)y^2-(a-2)(n+1)y-(a-2)=0\implies
$$}

\bigskip\noindent
{\bf Acknowledgements.} 
This research is funded by a Galatasaray University research grant and the T\"UB\.ITAK grants 115F412 and 113R017.

\bibliographystyle{amsplain}
\bibliography{references3}

 \end{document}